\documentclass[10pt]{amsart}
\usepackage{amsthm, amsmath, amssymb}
\usepackage[latin1]{inputenc}
\usepackage{graphicx}
\usepackage{amssymb}
\usepackage{amsmath}

\newtheorem{teo}{Theorem}
\newtheorem{Lemma}{Lemma}

\newtheorem{prop}{Proposition}

\newcommand{\m}{\mbox}

\theoremstyle{definition}
    \newtheorem{defin}{Definition}
    \newtheorem{es}{Example}
    \newtheorem{remark}{Remark}
\newtheorem{Notation}{Notation}

\begin{document}
\title{Algebraic functions with even monodromy}
\author{Michela Artebani}
\address{Dipartimento di Matematica ``DIMA'', Universit\'{a} di Genova, 
  via Dodecaneso 35,\linebreak 16146 Genova, Italia.}
\email{artebani@dimat.unipv.it}
\author{Gian Pietro Pirola}
\address{Dipartimento di Matematica ``F. Casorati'', Universit\'{a} di Pavia,
  via Ferrata 1,\linebreak 27100 Pavia, Italia.}
\email{pirola@dimat.unipv.it}
\date{\today}
\thanks{Partially supported by: 1) PRIN 2001:
  {\em Spazi di moduli e teoria di Lie};
  2) Gnsaga; 3) Far 2002 \\ \indent (PV):
  {\em Variet\'{a} algebriche, calcolo algebrico, grafi orientati e 
topologici.}}

\begin{abstract}
 Let $X$ be a compact Riemann
surface of genus $g$
and $d\geq 12g+4$ be an integer.
We show that $X$ admits meromorphic
functions with
monodromy group equal to the
alternating group $A_d.$
\end{abstract}
\maketitle


\section{Introduction}

Let $X$ be a compact Riemann surface of genus $g$ and $f:X \rightarrow {\mathbb{P}}^1$ be a meromorphic function of degree $d$.
The function field $C(X)$ is a finite algebraic extension of degree $d$ of $C({\mathbb{P}}^1)=\mathbb C(x)$.
The \emph{monodromy group} $M(f)$ is the Galois group associated to the Galois closure of the extension
$C(X)/\mathbb C(x)$.
The group $M(f)$ has a natural transitive representation in the symmetric group $S_d$.
In this paper, we prove the existence of meromorphic functions on $X$ with \emph{even monodromy}. This means that the monodromy group is contained in the \emph{alternating group} $A_d$ on $d$ elements. In fact, we study the case:
$$M(f)=A_d.$$
 
The problem of finding Riemann surfaces
with given monodromy group
is a classical one (see for example \cite{fried}).
In particular, a difficult question is that of determine all possible monodromy groups for the generic Riemann surface of genus $g$.
Several aspects of this problem were considered by Zariski in \cite{Z}. Notice that the definition of monodromy group can be extended to any holomorphic map between compact Riemann surfaces. In fact, Zariski  observed  that there are no holomorphic maps from the generic Riemann surface of genus $g>1$ to one of positive genus. Thus the critical case to study is that of meromorphic functions. Besides, he proved that the monodromy groups of a generic Riemann surface of genus $g>6$ are not solvable. Zariski also introduced the reduction to the \emph{primitive} case. This means that we can consider only minimal extensions of $\mathbb C(x)$ or, equivalently, \emph{indecomposable} meromorphic functions.

It was first shown by Guralnick and Thompson in \cite{GT} that there are groups that cannot occur as the monodromy group of a meromorphic function on a Riemann surface of genus $g$.

Later, it was proved that the monodromy group of an indecomposable meromorphic map of degree $d$ on a generic Riemann surface of genus $g>3$ is either $S_d,$ $d\geq (g+2)/2$ or $A_d,$ $d\geq 2g+1$.
This result is contained in a series of three papers: \cite{gu2}, \cite{GM} and \cite{GS}.

It is well known that the symmetric group of degree $d\geq (g+2)/2$ can be the monodromy group of any curve of genus $g$ (see \cite{KL}).
As noticed by Guralnick in \cite{gu}, this makes the case of the alternating group very interesting.

In a recent article
K. Magaard 
and H. V\"olklein prove that the general curve of genus $g\geq 3$ admits meromorphic functions with monodromy group $A_d$ and 3-cycles (or double transpositions) as branch cycles 
if and only if 
$d\geq 2g+1$ (see \cite{mv}). 

In this paper we prove that, for $d$ sufficiently large compared to $g$, every compact Riemann surface of genus $g$ realizes the monodromy group $A_d.$

\begin{teo}\label{teo}
Let $X$ be a
compact Riemann surface of genus $g>0$
and $d\geq 12g+4$ be an integer.
Then there exist indecomposable meromorphic functions $f\in \mathbb C(X)$ of degree
$d$
with $M(f)=A_d.$
Moreover,
there exists on $X$ a family
$\mathcal F(X,d)$ of such maps with dimension:
$$
\dim\mathcal F(X,d)=[\frac{d+3}{2}]-2g+2.
$$
\end{teo}

We give a brief summary of the proof,
which is contained in section $4.$

We fix  a \emph{spin bundle} $S$
on a compact Riemann surface $X.$
This is  a line
bundle on $X$ such that $S^2=K_X,$
where $K_X$ is the canonical bundle.
We then take a divisor $D$ on $X$
   having large degree and small support
(i.e. $D=n_1P_1+n_2P_2+n_3P_3$ where $P_1,$$P_2$
   and $P_3$ are points of $X$).
Let $[D]$ be the divisor given by the formal sum of points in the support of $D.$
The square of a global section
$s$ of $S(D)$ is a meromorphic form
   $\omega=s^2$ on $X$
with at most poles at $[D].$
It has even order at each point.
If $\omega$ is \»exact" i.e.
$\omega=df,$ then $f$ is a
meromorphic function on $X$
with local monodromy in
   $A_d$ where $d=\deg (f).$
We call such meromorphic functions  \emph {odd ramification
   coverings} of the projective line.
\footnote{
We avoid the tempting term
``spin covering'' from
a suggestion of M. Fried who explained
    us a different use of this word.}
These coverings were studied in \cite {fried3}, \cite {se 1}
and \cite {se 2} in connection with the
spin lifting.

Given an inferior bound on the degree of $D$, the existence of not trivial sections $s$ of $S(D)$ with
exact square follows easily for dimensional reasons.
In section $4,$ we compute the dimension
   of the families
of these sections when $d\geq 12g+4.$
As in the theory of special
divisors (see \cite{cor}),
it is possible
to compute the tangent
   space to our families in the
cohomological setting.
For this, we use a
rational variation of
De Rham's algebraic theory
   and the adjustment of a
technical result contained in \cite{pi}
(we recall its proof 
in section $6$).
The generic smoothness of our
families follows from a
cohomology vanishing.
Thus, comparing dimensions,
we find indecomposable functions with
the expected degree:
their poles have maximum order at $[D].$
An algebraic lemma and a count of parameters show that the
monodromy group of
the general element
in the family is exactly $A_d.$

In section $5,$ using the theory
of admissible coverings,
we show that odd ramification coverings are
limits of \emph{simple
odd ramification coverings}
   i.e. with ramification
   points all of order $3$ with distinct
images in $\mathbb P^1.$
Then, Theorem \ref{teo} implies that
the generic Riemann surface of genus $g$  admits  simple
odd ramification coverings
   of degree $d\geq 12g+4$
with monodromy group $A_d$ (see Theorem 2).
This proves in a completely different way a weaker version of the result contained in \cite{mv}.\\

In \cite{bp} S. Brivio and the second author apply Theorem \ref{teo} to give an analogous result for rational functions on complex surfaces. \\

\noindent\textbf{Acknowledgement}. We would like to thank the referee and Michael D.
Fried for valuable suggestions. We are grateful to Juan Carlos Naranjo and Enrico Schlesinger for many helpful conversations.
 




\section{Preliminaries and Notation}

We start recalling some known results about the monodromy group of a meromorphic function
(see \cite{mi}).
Let $X$ be a compact Riemann
surface of genus $g,$
$f:X\rightarrow \mathbb P^1$
   be a meromorphic function of
degree $d$
on $X$ with branch locus $B$ and ramification locus 
$R$. Let $p:X\setminus R \rightarrow \mathbb P^1\setminus B=Y$
be the associated topological covering.
Fixed a base point $y\in Y$, the fundamental group $\Pi_1(Y,y)$ acts on the fibre $f^{-1}(y)$ 
by path lifting giving a transitive subgroup of the symmetric group $S_d$.
This is the \emph{monodromy group} of $f$ (it is determined by $f$ up to conjugacy).
In fact, the monodromy group of $f$ is isomorphic to the Galois group associated to the Galois closure of the extension $\mathbb C(X)/\mathbb C(x)$ (see \cite{h}).

Conversely, let us consider $B=\{b_1,\dots,b_r\}\subset \mathbb P^1$, 
$Y=\mathbb P^1\setminus B$ and fix
standard geometric generators of $\Pi_1(Y,y)$. Riemann's existence theorem gives a bijection between 
meromorphic functions
of degree $d$ with branch locus $B$ (up to isomorphism) and ordered r-uples of
permutations $\{\sigma_1 ,\dots,\sigma_r\}$ with $\sigma_1\cdot\cdot\cdot\sigma_r=id$
and generating a transitive
subgroup of $S_{d}$ (up to conjugacy).

It can be easily seen that $f$ is indecomposable if and only if its 
monodromy group is primitive.

We need the following classical result:  
\begin{Lemma}\label{l1}
Let $G$ be a transitive and
primitive subgroup of the alternating group $A_d$ containing a
3-cycle.
Then: $G=A_d.$
\end{Lemma}
These examples provide a simple way to state indecomposability:
\begin{es} \label {es} \ \\
A)\ Let $f:X \rightarrow \mathbb P^1$ be a
decomposable meromorphic map of degree $d$ and $P,Q\in X$ with
$\{P,Q\}=f^{-1}(\infty)$ and ramification
indices $e_P,e_Q>1$.
Then $e_P$ and $e_Q$ are not relatively prime.
\begin{proof}
Let $f=f_2\cdot f_1$ be a factorization with $f_1:X\rightarrow Y$ and $f_2:Y\rightarrow \mathbb P^1$ of degrees $d_1,d_2>1$. The fiber of $f_2$ over $\infty$ consists at most of two points.

If ${f_2}^{-1}(\infty)=\{P',Q'\}$ then $P$ and $Q$ are total ramification points for $f_1$. 
Thus $e_P=d_1p'$ and $e_Q=d_1(d_2-q')$, where $p'$ is the ramification index of $f_2$ in $P'$.

If ${f_2}^{-1}(\infty)=\{R'\}$ then $R'$ is a total ramification point for $f_2$. This implies that $e_P=pd_2$ and $e_Q=(d_1-p)d_2$ where $p$ is the ramification index of $f_1$ in $P$.
\end{proof}
\noindent B)\ Let $f:X \rightarrow \mathbb P^1$ be a decomposable
meromorphic map of degree $d$ and $P,Q,R\in X$ with
$\{P,Q,R\}=f^{-1}(\infty)$ and odd ramification
indices $e_P,e_Q, e_R>1$.
Then $e_P$, $e_Q$ and $e_R$ are not relatively prime.
\begin{proof}
There are three cases according to the cardinality of $f_2^{-1}(\infty)$. Two of them are analogous to example A. We only consider the case $f^{-1}(\infty)=\{P',Q'\}$. 

Let $p'$ be the ramification index of $f_2$ in $P'$. 
A unique point in $f^{-1}(\infty)$ is totally ramified for $f_1$. Suppose that $R$ is such a point and that $f_1(R)=Q'$. Let $p,q$ be the ramification indices of $f_1$ in $P,Q$. 
Then we have: $e_P=pp'$, $e_Q=qp'$ and $e_R=d_1(d_2-p')$.
Notice that $d=d_1d_2=e_P+e_Q+e_R$ is odd, in particular $d_2$ is odd. Besides, $p'$ is odd because $e_P$ is odd. This contradicts the fact that $e_R$ is odd too, then this case can not occur.
\end{proof}
 \end{es}
\begin{Notation}
Let $D$ be a divisor on $X$, we denote by $\mbox{supp}(D)$ its support and by $[D]$ the divisor given by the formal sum of points in $\mbox{supp}(D)$. Let $X\backslash [D]$ be the open set $X\backslash \mbox{supp}(D)$. 
We write $[D_1]\cap[D_2]=\emptyset$ if $D_1$ and $D_2$ have disjoint supports.
If $D_1-D_2$ is effective we use the notation $D_1-D_2\geq 0$.
\end{Notation}




\section{Spin sections with exact square}

\begin{defin}
Let $X$ be a compact Riemann surface,
we call \emph{spin bundle}  a line
bundle $S$ on $X$ with:
$S^2=K_X,$
where $K_X$ is the canonical bundle.
\end{defin}
It is known (see \cite{cor}) that there are exactly $2^{2g}$ non
equivalent spin bundles on any compact Riemann
surface $X$ of genus $g$.

Fix $X$ a compact Riemann surface of
genus $g,$ $S$ a
spin bundle on $X$ and $D\in Div(X)$
a divisor of the form:
$$D=n_1P_1+n_2P_2+n_3P_3\ \ \ (\star) $$
where $n_1,n_2,n_3$ are integers with $n_1\geq n_2\geq n_3\geq 0$ and
$P_1,P_2,P_3$ are distinct points in $X$.
Let $k$ be the degree
of $[D]$. Notice that $1\leq k\leq 3$ and $D=\sum_{i=1}^{k}n_iP_i$.
We set: $$D'=2D-[D].$$
Consider now the line bundle $S(D).$
The squares of sections
$s\in H^0(X,S(D))$ give
meromorphic forms
$\omega=s^{2} $
on $X$ with at most poles in $[D].$
More precisely, their divisor is:
$$(\omega)_{0}-(\omega)_{\infty}=2div(s)-2D=2E,$$
with $E=div(s)-D$. 
If $\omega$ is an exact form with $\omega=df$ then
$f\in H^{0}(X,\mathcal O_X(D'))$ and
\begin{eqnarray}
e_f(x)&=&ord_{\omega}(x)+1
=2ord_E(x)+1\ \ \m{if}\ x\in (\omega)_{0}\nonumber,\\
e_f(x)&=&-ord_{\omega}(x)-1=
-2ord_E(x)-1\ \ \m{if}\ x\in (\omega)_{\infty}
     \nonumber,\end{eqnarray}
where $e_f(x)$ is the ramification index of $f$ in $x$. 
The map $f$ has odd ramification
index in each point, hence
\emph{even monodromy},
i.e. $M(f)$ is contained
in the alternating group.
In fact  each permutation
generating the
monodromy group $M(f)$
can be decomposed
in cycles of odd length.
\begin {defin} We call a meromorphic map
   $f:X\rightarrow \mathbb P^1$ an
   \emph{odd
ramification covering}
    if all ramification
    points have odd index.
\end{defin}

We define the analytic map:
$$
\Upsilon: H^0(X,S(D))\rightarrow
   H^1(X\backslash [D],\mathbb C)\ \ \
s\mapsto [s^2].
$$

Firstly, we want to show
the existence of
sections $s\in H^0(X,S(D)),s\not= 0$
with $ \Upsilon(s)=0.$
\begin{Notation}
   $h^i(X,L)=\dim H^i(X,L)$ (i=1,\ 2);\
$\mathcal H(X,D)=\Upsilon^{-1}(0)$; \\$\mathcal F(X,D)=\{f\in\mathbb
C(X):df=s^2,\ s\in \mathcal H(X,D)\}$.
\end{Notation}
By the Riemann-Roch theorem and
De Rham theory we have :
$$
h^0(X,S(D))=\deg(D)\
{ \rm if}\  \deg(D)>g-1;$$
$$ h^1(X\backslash[D],\mathbb C)=2g+k-1
$$
where $k$ is the degree of $[D]$.

Suppose $\Upsilon(s)\not=0$ for every
$s \in H^0(X,S(D))\setminus \{0\}.$
This gives a contradiction
   if $\deg(D)>2g+k-1.$
Therefore we have:

\begin{prop}
\label{dim}
Let $X$ be a compact Riemann surface
   of genus $g$ and $D$ a divisor
as in $(\star).$
If $\deg(D)>2g+k-1\ $ then
$\dim_{\mathbb C}\mathcal H(X,D)>0.$
In particular, there exists a not constant
odd ramification covering $f\in \mathcal F(X,D)$.
\end{prop}

To get more precise information
about the dimension of $\mathcal H(X,D),$
we are interested in the surjectivity
of the differential of $\Upsilon$
in $s\in \mathcal H(X,D)$:
$$
d\Upsilon(s):H^0(X,S(D))\rightarrow
H^1(X\backslash[D],\mathbb C).
$$
Take the holomorphic curve $s(t)=s+tv$ in $H^0(X,S(D))$ then:
$$
[s(t)^2]=[s^2+2tsv+t^2v^2]=2t[sv]+o(t^2).
$$
Therefore:
$$
   d\Upsilon(s)(v)=[2sv].
$$

We now recall a technical lemma
   contained in \cite{pi},
its proof is contained 
in section $6.$
Let $A,B\in Div(X)$
be effective divisors with $[A]\cap[B]=\emptyset$, 
$a=\deg[A]$ and $b=\deg[B]$.
Consider the map:
$$
\mathcal D: H^0(X,\Omega(A-B))\rightarrow
H^1(X\setminus [A],\mathbb
C)\ \ \ \omega\mapsto
[\omega].
$$
\begin{Lemma}\label{tec}The map $\mathcal D$ is onto if
$\deg(A-B)-a-b>2g-2$.
\end{Lemma}
We consider a divisor $D$
as in $(\star)$ and a section
$s \in \mathcal H(X,D)\backslash \{0\}$
with zero divisor $div(s)=E.$
The multiplication by $s$
gives an isomorphism:
$$
m(s):H^0(X,S(D))\rightarrow H^0(X,\Omega(2D-E)).
$$
Composing  with $\mathcal D$ 
we obtain:
$$
\mathcal D\circ m(s):
H^0(X,S(D))\rightarrow H^1(X\setminus
[D]),\mathbb C)\ \ \
v\mapsto [sv],
$$
which is
$d\Upsilon(s)$ (up to a costant).
\smallskip

We would like to apply Lemma \ref{tec}.
We first decompose the
divisor of $s$ in the sum of two effective divisors as follows:
$$
E=E_1+E_2\ \m{with}\
[D]-[E_1]\geq 0,\ [E_2]\cap [D]=\emptyset.
$$
Let us define:
$$A=2D-E_1,\ B=E_2.
$$

With this choice, Lemma \ref{tec}
gives the surjectivity of $d\Upsilon(s)$
if the couple $(D,E)$ satisfies :

\begin{enumerate}
\label{pr}
\item[a)]
$2D-E_1\geq 0,$
\item[b)]$\deg(2D-E)-k-\deg[E_2]>2g-2$
\end{enumerate}
\begin{Notation}
Let $T \in X\setminus[D]$ be a generic point and $r\geq 1$
be an integer. Set
\begin{eqnarray}
\mathbb P(\mathcal H(X,D)):=
\{ (s)\in \mathbb PH^0(X,S(D))\ :
   s\in \mathcal H(X,D)\}
   ,\nonumber\\
 \ \ \ \ \ \mathbb PH^0(X,S(D-rT)):=\{ (s)\in
\mathbb PH^0(X,S(D))\ : div(s)\geq rT\}.
\nonumber
\end{eqnarray}
\end{Notation}
Notice that:
\begin{eqnarray}
&\ &\dim\mathbb P(\mathcal H(X,D))\geq
   \deg(D)-2g-k,\label{dim1}\nonumber\\
&\ &\dim \mathbb PH^0(X,S(D-rT))=\deg
(D)-1-r\label{dim2}
.\nonumber\end{eqnarray}

Therefore, if $r\leq \deg(D)-2g-k:$
$$
\mathbb P(\mathcal H(X,D))\cap\mathbb P
H^0(X,S(D-rT))\not=\emptyset.
$$
Let $s\in\mathcal H(X,D)$
be a not trivial section with
$div(s)\geq(\deg (D)-2g-k)T.$
   Then we have:
$$
\deg(E_1)\leq (g-1+\deg(D))-(\deg(D)-2g-k)=3g+k-1.
$$
In particular a)
holds for $E=div(s)$ when:
$$2n_i>3g+k-1,\ i\in\{1,\dots,k\}.$$
Then  $b)$ becomes:
$\deg(D)-\deg[E_2]>3g-3+k.$
Since
$$\deg[E_2]\leq \deg[E]\leq 3g+k-1,$$
condition b) holds when
$$\deg(D)>6g+2k-4.
$$
Under these hypotheses, Lemma \ref{tec}
implies that all components of
$\mathcal H(X,D)$ contain
   a point
$s$ such that $d\Upsilon(s)$ is
surjective.
The implicit function
theorem finally gives:

\begin{prop}\label{su}
Let $X$ be a compact Riemann surface
of genus $g$ and
$D$ be a divisor as in $(\star)$
   with support of degree
$k.$ Suppose that:
\begin{itemize}
\item[1)]$2n_i>3g+k-1,\  i\in\{1,\dots,k\}$
\item[2)]$\deg(D)>6g+2k-4.$
\end{itemize}
Then the general point of any component of
$\mathcal H(X,D)$ is smooth and
$$
\dim \mathcal H(X,D)=\deg(D)-2g-k+1.
$$
\end{prop}




\section{Constructing maps with even monodromy}

In the previous section we have seen that
sections in  $\mathcal H(X,D)$
give in a natural way odd ramification coverings.
However Proposition \ref{su} only gives
a bound on the ramification at
infinity of such maps. In fact ``a priori"
these sections could have zeros at $D$,
  hence the degree of the corresponding maps could drop.
Moreover,
we need to provide indecomposable maps.
\bigskip

Let $D=n_1P_1+n_2P_2+n_3P_3$ as in the previous sections with $k$ the cardinality of its support. Set:
$$d_1=2n_1-1,\ d_2=2n_2-1,\ d_3=2n_3-1.$$

\begin{defin} We will say that $(d_1,d_2,d_3)$
is an \emph{indecomposable triple} in these cases:
\begin{itemize}
\item if $k=1$: $d_1$ is prime number;
\item if $k=2$: $d_1,d_2$ are relatively prime,
\item if $k=3$: $d_1,d_2,d_3$ are relatively prime.
\end{itemize}
\end{defin}

\noindent The result is the following:

\begin{prop}\label{fin}
Let $X$ be a compact Riemann surface of
   genus $g>0$ and
$D$ be a divisor as in $(\star)$
with support of degree $1\leq k\leq 3.$
   Let $d_1,d_2$ and
    $d_3$ as before. Suppose:
    \begin{enumerate}
\item[a)] $d_i>3g+k$ for all $i\in\{1,\dots, k\}$
   \item[b)]$\deg(D)>6g+2k-3,$
    \item[c)]$(d_1,d_2,d_3)$ is an indecomposable triple.
      \end{enumerate}
Then there exists a family $\mathcal F(X,D)$ of meromorphic functions $f$
on $X$ with these properties:
\begin{itemize}
\item[1)] $deg(f) = 2\deg(D) -k,$
\item[2)] $M(f)\subset A_d$ and $M(f) $ is primitive i.e. $f$ is
indecomposable.
\item[3)] $f$ has maximum ramification in $[D].$
\end{itemize}
Moreover:   $\dim\mathcal F(X,D)=\deg(D)-2g-k+2.$
\end{prop}

\begin{proof}

By construction the divisors
$$D_i=D-P_i,\ i\in\{1,\dots,k\}$$
satisfy the hypothesis of
   Proposition \ref{su}.
   Then we get:
$$
\dim \mathcal H(X,D_i)=\deg(D)-2g-k=
\dim\mathcal H(X,D)-1,\ \ i\in\{1,\dots,k\}.
$$
Comparing dimensions,
we get the existence of sections $s$ in:
$$ \mathcal H(X,D)
\setminus\mathcal \bigcup_{i=1}^{k} \mathcal
H(X,D_i)
.$$
They have the property:
$$
\ ord_{P_i}(s^2)=-2n_i,\ i\in\{1,\dots,k\}.
$$

This proves the existence of a family of
odd ramification coverings $\mathcal F(X,D)$ with
maximum ramification in the points of $[D].$
Notice that maps in $\mathcal F(X,D)$
   have degree: $$d=2\deg(D)-k.$$
Moreover:
$$\dim\mathcal F(X,D)=\deg(D)-2g-k+2.$$

Finally, condition c) and a previous
   observation (see section 2, Example \ref{es})
give indecomposability of maps in $\mathcal F(X,D).$
\end{proof}

Proposition \ref{fin} specializes to three cases as follows:
\begin{itemize}
\item $k=1$: we get families of indecomposable
odd ramification coverings
with a pole $P_1$ of total
ramification and prime degree $d.$
\item $k=2$: we get families of indecomposable
   odd ramification coverings of (even) degree
   $d=2n_1+2n_2-2$ and ramification $d_{1}=2n_1-1$
   at $P_1,$ $d_{2}=2n_2-1$ at $P_2$
with $d_{1}$ and $d_{2}$ relatively prime.
\item $k=3$:   we get families of indecomposable
   odd ramification coverings of (odd) degree
   $d=2n_1+2n_2+2n_3-3$ and ramification $d_{1}=2n_1-1$ at $P_1,$
    $d_{2}=2n_2-1$ at $P_2,$
   $d_{3}=2n_3-1$ at $P_3$ with $d_{1},d_{2},d_{3}$ relatively prime.

\end{itemize}
Notice that condition $b)$ in Proposition \ref{fin} implies:
$$d>12g+3k-6.$$
It can be easily proved that, for every integer $d\geq 12g+4$, there exists a divisor $D$
in the hypotheses of Proposition \ref{fin} with $d=2\deg(D)-k$.
Fixed $d$ a positive integer, we use the following notation:
$$\mathcal F(X,d)=\bigcup_{D\in\mathcal D(d)}
\mathcal F(X,D),$$
where $\mathcal D(d)=\{D\in Div(X): D\ {\rm as\ in\ Proposition\ \ref{fin}},\
k>1,\  2\deg(D)-k=d\}$.
Notice that, for every $D\in \mathcal D(d)$:
$$\dim\mathcal F(X,d)=\deg(D)-2g+2=[\frac{d+3}{2}]-2g+2.$$
\begin{proof}[Proof of Theorem 1]
Notice that Proposition \ref{fin} gives a
weak version of Theorem 1.
In fact, we need to show that the generic
map in $\mathcal F(X,d)$ has monodromy group equal to $A_d.$
By Lemma \ref{l1},
it is enough to prove that its monodromy
group contains at least
a 3-cycle.
Let $d$ be an even integer.
   Suppose, by contradiction,
   that the monodromy group
   of a generic map $f\in\mathcal F(X,d)$
   doesn't contain any $3$-cycle. It follows that any branch
   point has ramification 
   order $\geq 4.$
Then Hurwitz's formula gives:
$$
2g-2=-2d+(d-2)+\sum_{P\in X}(e_f(P)-1)
\geq -d-6+4b
$$
where $b$ is the number of branch points.
Thus:
$ b\leq (2g+4+d)/4.$
Then, from the theory of Hurwitz schemes:
$$\dim\mathcal F(X,d)\leq (2g-4+d)/4.$$
On the other hand:
$$\dim\mathcal F(X,d)=d/2-2g+2>(2g-4+d)/4\
\m{ if }\ d>10g-12.$$
Similar considerations hold for odd $d.$
\end{proof}

\begin{remark}
Let $\mathcal F_1(d)$ be the union of all families $\mathcal F(X,D)$
with $k=1,$ $2\deg(D)-1=d$ and $X$
moving in the moduli space of smooth projective curves $\mathcal M_g.$
A parameters' count shows that the generic
element in $\mathcal F_1(d)$
is a couple $(X,f)$ where $f$ has total
ramification at infinity and all
other ramification points of index $3$ with
distinct images in $\mathbb P^1.$
\end{remark}




\section{Simple odd ramifications}

Let $X$ be a compact Riemann
surface of genus $g$ and
    $f:X \to \mathbb P^{1}$ a meromorphic function
     with branch
locus $B$ and ramification locus $R.$
The simplest odd ramification type
is a $3:1$ ramification covering
(see \cite{fried3}).
\begin {defin} We say that an
odd ramification covering is
\emph{simple} if all points in $R$
have ramification index equal to $3$ and
   have distinct images in
$\mathbb P^{1}$
\end{defin}
Using admissible coverings we prove:
\begin{prop}\label{defor}
1) Any  odd ramification covering $f:X\rightarrow \mathbb P^1$
is the specialization of a
family of simple odd ramification coverings
$f_{s}:X_s\rightarrow \mathbb P^1.$\\2)
The monodromy group of $f$ is
   a subgroup of the monodromy group of
$f_{s}$ for general $s.$
\end {prop}
\begin{Notation}Set $Y_{0}= \mathbb P^{1}$ and let
$T\in Y_{0}$ be a branch point of $f.$
Let $Y_1$ be another copy of
$\mathbb P^{1}$ and $P$
a point in $Y_{1}.$
We call  $\tilde Y$
   the singular rational curve obtained by
   gluing  $Y_{1}$ and $Y_{0}$
   in $T$ and $P$:
     $\tilde Y= Y_{0}\cup Y_{1}.$
Let $D= f^{-1}(T)= \sum_{i=1}^{r}(2n_{i}-1 )P_{i}.$
For any
$i\in\{1,\dots,r\}$ we take
a copy $X_{i}$ of $\mathbb P^1$
with a selected point
$Q_{i}\in X_{i}.$
Let
$$
\tilde X=(\cup_{i=1}^{r} X_{i})\cup X
$$
be the curve obtained
identifying $P_i$ and $Q_i$ for every
$i\in\{1,\dots,r\}.$
\end{Notation}
\begin{Lemma}
We can define a map 
$
\tilde f : \tilde X \rightarrow \tilde Y
$
satisfying the following properties:

\begin{enumerate}
\item[a)] $\tilde f|_{X} =f$ ,
   \item[b)] $\tilde f|_{X_{i}}=f_i:X_{i}\to Y_{1}$
for every $i\in\{1,\dots,r\},$
\item[c)] $f_{i}$ is an odd ramification covering of
degree $2n_{i}-1$ with total
ramification at $Q_{i},$
\item[d)] ramification points of $f_i$
other than $Q_i$ have index $3,$
\item[e)] branch points of the $f_{i}$
other than  $T$ are all distinct.
\end{enumerate}
\end{Lemma}
\begin{proof} The existence
of meromorphic functions
$X_{i}\rightarrow \mathbb P^1$
satisfying c),d) and e) follows from
Riemann's existence theorem.
\end{proof} 
Notice that we have constructed an admissible covering.
Now we can apply the standard
smoothing procedure (see \cite{hm} and \cite{cor}).

\begin{proof}[Proof of Proposition \ref{defor}]
We can consider $\tilde Y$ as the union of two
lines in the projective plane: $\tilde Y=\{xy=0\}.$
Choose local coordinates
$w_{i}$ and $z_{i}$ at $Q_{i}=P_{i}$
such that locally $\tilde X$
is given by $\{w_{i}z_{i}=0\}$ and
$\tilde f$ has the form
$x= w_{i}^{2n_{i}-1}$ on $X$ and
$y=z_{i}^{2n_{i}-1}$ on $X_{i}$
for each $i\in\{1,\dots,r\}.$
Set now: $m_{i}=2n_{i}-1$ for $i=1,\dots,r$. 
Let $\Delta$ be the curve in $\mathbb C^{r+1}$
defined by:
$$
\Delta=\{(t,t_1,\dots,t_r):t=t_{i}^{m_i},
\mid t\mid<1,\mid t_{i} \mid <1, i=1,\dots,r\}.
$$
For every
$s=(t,t_1,\dots,t_r)\in \Delta$
consider:
$$
I_s=\{(x,0):\mid x\mid\leq \mid t(s)\mid\}\cup\{(0,y):
\mid y\mid\leq \mid t(s)\mid\}\subset \tilde Y,
$$
$$
J_{i,s}=\{(w_i,0):\mid w_{i}\mid\leq \mid t_{i}(s)\mid\}
\cup\{(0,z_{i},s):
\mid z_{i}\mid\leq \mid t_{i}(s)\mid\}\subset \tilde X.
$$
Let $U_{s}=\tilde Y\setminus I_s$
and $U'_s=\tilde X\setminus
\cup_{i} J_{i,s}.$
Finally, we define:
$$
V_s=\{(x,y)\in \mathbb C^2: xy=
t(s),\mid x\mid<1,\mid y\mid<1\},
$$
$$V_{i,s}=\{(w_{i},z_{i})\in
\mathbb C^2: w_{i}z_{i}=
t_{i}(s),\mid w_{i}\mid<1,\mid z_{i}\mid<1\}.
$$
We smooth locally $\tilde Y$ gluing
$U_s$ to $V_s$ in the following way:
$$
(x,0)\mapsto (x,\frac{t(s)}{x});
\ (0,y)\mapsto (\frac{t(s)}{y},y).
$$
We can do the same for $\tilde X$ gluing $U'_s$ to $V_{i,s}$:
$$
(w_i,0)\mapsto (w_i,\frac{t_i(s)}{w_i});
\ (0,z_i)\mapsto
(\frac{t_i(s)}{z_i},z_i).
$$
For each $s\in \Delta$ we get
a compact smooth surface $X_s$ and a copy $Y_s$ of
$\mathbb P^1$. We can define a map
$f_{s}:X_{s}\to Y_s$ with
$f_s=\tilde f$ on $U'_s $ and
$f_{s}|V_{i,s}$ given by:
$$
(w_i,z_i)\mapsto(w_{i}^{m_i},z_{i}^{m_i}).
$$

In fact, we get compact
Riemann surfaces $X_{s}$
(the complex structure on $X_{s}$
is given by the ramified covering)
of genus $g$ and odd ramification coverings
$f_{s}$ on $X_{s}$
   with the same degree of $f$
that specialize to $ \tilde{f}$ for
$s=0.$
Remark that the monodromy group of $f$
is determined by the branches in
$Y_0\setminus \{T\}.$ Let
$\gamma_1,\dots,\gamma_b$ be loops in
$Y_0$ corresponding to these branches.
Since the coverings
$f_{s}$ and $f$ are topologically
the same in $U'_s,$ we can find
$\gamma_{1,s},\dots,\gamma_{b,s}$
loops in $Y_s\setminus V_s$ giving
the same monodromy as
$\gamma_1,\dots,\gamma_b.$
Then, the monodromy group $M(f_{s})$ contains
$M(f)$ as a subgroup for every
$s\in\Delta.$ Moreover,
the coverings $f_{s}$ have simpler odd
ramification.
Repeating this procedure for each branch point 
we finally get a family of simple
odd ramification coverings deforming $f.$
\end{proof}

\begin{remark}
It is a not trivial problem to find a deformation of $f$ as in Proposition \ref{defor} which preserves the conformal structure of $X.$
\end{remark}
Combined with Theorem 1 this implies:
\begin{teo}
The generic Riemann surface of genus $g$
admits simple odd ramification coverings of
degree $d\geq 12g+4$ with monodromy
group equal to $A_d.$
\end{teo}



\section{A technical result}

Let
$A,B\in Div(X)$ be two effective
divisors on $X$ with $[A] \cap [B]=\emptyset,$
$a=\deg[A],$ $b=\deg[B].$
Let us call $A'=A-[A].$
We define a sheaf
$\mathcal R=\mathcal R(A,B)$
as follows:
$$\mathcal R(A,B)(U):=\{f \in
\mathcal O(A')(U):\ df\in \Omega(-B)(U)\},
$$
$U\subset X$
open set.
We have the short exact
sequence of sheaves:

\begin{equation}
\label{1}0\rightarrow \mathcal
O(A'-B-[B])\stackrel{i}{\rightarrow} \mathcal R
\stackrel{e_v}{\rightarrow} \mathbb C_{B}
\rightarrow 0
\end{equation}
where $\mathbb C_{B}$ is the
skyscraper sheaf with support in
$[B]$ and $e_v$ is the
\»valutation'' morphism.
If $f \in \mathcal R(A,B)(U)$
then $df \in \Omega(A-B)(U),$
a second exact sequence is then:
\begin{equation}
\label{2}
0\rightarrow \mathbb C
\stackrel{i}{\rightarrow} \mathcal
R\stackrel{d}{\rightarrow}
\Omega(A-B)^{\circ}\rightarrow 0
\end{equation}
where $\Omega(A-B)^{\circ}$
is the subsheaf of forms with residues
zero in $[A].$
Finally, taking residues,
we get a third sequence:
\begin{equation}
\label{3}
0\rightarrow\Omega(A-B)^{\circ}\rightarrow
    \Omega(A-B)\stackrel{res}{\rightarrow}\mathbb C_{A}
   \rightarrow 0.
   \end{equation}
 

Consider the coboundary
operators associated to sequence $(2)$:
$$\Delta_i:H^{i-1}(X,\Omega(A-B)^{\circ})\rightarrow H^i(X,\mathbb
C)\ \ i=1,2.$$
Suppose now:
$H^1(X,\mathcal R)=0.$
Then $\Delta_1$ is surjective and $\Delta_2$ is an isomorphism.
Moreover, the image of the  map
$res:H^0(X,\Omega(A-B))
\rightarrow H^0(X,\mathbb C_A)=\mathbb C^a$
is contained in $\mathcal I=
\{\sum_{1}^{a}x_{i}=0\}$ by
the residue theorem.
Sequence $(3)$
gives then: $H^1(X,\Omega(A-B))=0.$
Therefore $\Delta_1 $ and
$res$ are both surjective. Finally,
this implies the surjectivity of
the De Rham's map:
$$
\mathcal D: H^0(X,\Omega(A-B))
\rightarrow H^1(X\setminus [A],\mathbb
C)\ \ \ \omega\mapsto
[\omega].
$$

\begin{prop}
\label{R}
If $H^1(X,\mathcal R)=0,$ then
the map $\mathcal D$ is surjective.
\end{prop}
In fact, the vanishing of
$H^1(X,\mathcal R)$ and the surjectivity of
$\mathcal D$
    are equivalent conditions (see \cite{pi}).

\begin{proof}[Proof of Lemma \ref{tec}]
We have $H^1(X,\mathbb C_{B})=0$ and
$H^1(X,\mathcal O(A'-B-[B])=0.$ This implies
$H^1(X,\mathcal R)=0.$

\end{proof}




                                   %
                                   %
                                   %
                                   %
                                   %
                                   %

\end{document}